\documentclass[12pt]{article}
\usepackage{amsmath}
\usepackage{amsfonts}
\usepackage{amssymb}
\newcommand{\BB}{{\bf B}}

 \bigskip\bigskip

\begin{document}
\author{Xiaojun Huang, Shanyu Ji and  Dekang Xu}
\date{}
\title{ A New Gap Phenomenon for  Proper Holomorphic Mappings   from ${\BB}^n$ into ${\BB}^{N}$}
\maketitle

In this paper (Math. Res. Lett. 13 (2006). No 4, 509-523),
  the authors established a pseudo-normal form for proper holomoprhic
mappings between balls in complex spaces with
degenerate rank.  This then was used to give a complete characterization
for all proper holomorphic maps with geometric rank one,
which, in particular, includes the following as an immediate application:

Theorem: Any rational holomorphic map from $B^n$ into $B^N$ with $4\le n\le
N\le 3n-4$ is equivalent to the D'Angelo map
$$F_{\theta}(z',w)=(z',(\cos\theta)w,(\sin\theta)z_1w , \cdots,
(\sin\theta)z_{n-1}w, (\sin\theta)w^2, 0'),\ 0\le  \theta\leq \pi/2.$$

It is a well-known (but also quite trivial) fact that any non-constant
rational CR map from  a piece of the sphere $\partial {B^n}$ into the
sphere $\partial {B^N}$ can be extended as a proper rational holomoprhic
map from $B^n$ into $B^N$ ($N\ge n\ge 2$).  By using the rationality
theorem that the authors
established in  [HJX05], one sees that the the above theorem (and also the
main theorem of the paper)   holds in the same way for any non-constant
$C^3$-smooth CR map from a piece of  $\partial {B^n}$ into
$\partial{B^N}$.

The   paper [Math. Res. Lett. 13 (2006). No 4, 509-523]
was first  electronically  published  by Mathematical
Research Letters   several months ago  at its home website:

http://www.mrlonline.org/mrl/0000-000-00/Huang-Ji-Xu2.pdf.

(The pdf file of the printed journal version  can also be downloaded at

http://www.math.uh.edu/~shanyuji/rank1.pdf).

 \bigskip

  \centerline{\bf References}

\noindent
[A77] H. Alexander,
Proper holomorphic maps in ${\bf C}^n$,  Indiana Univ. Math. Journal 26,
137-146 (1977).

\noindent
[BER99] M. S. Baouendi, P. Ebenfelt and L. Rothschild,
Real Submanifolds in Complex Spaces and Their Mappings,
Princeton Univ. Mathematics Series 47, Princeton University, New
Jersey, 1999.

\noindent
[BR90]  M. S. Baouendi,  and L. P.
Rothschild, Geometric properties of mappings between
hypersurfaces in complex spaces, J. Differential Geom. 31, 473-499, 1990

\noindent
[CS90] J.Cima and T. J. Suffridge, Boundary behavior of rational proper
maps, Duke Math. J. 60, 135-138 (1990).

\noindent
[DA88] J. P. D'Angelo, Proper holomorphic mappings between balls of
different dimensions, {\it Mich. Math. J.} 35, 83-90 (1988).

\noindent
[DA93] J. P. D'Angelo, Several Complex Variables and the Geometry of
Real Hypersurfaces, CRC Press, Boca Raton, 1993.

\noindent [CD96] J. D'Angelo and D. Catlin,
 A stabilization theorem for
Hermitian forms
and applications to holomorphic
mappings, Math Research Letters 3, 149-166 (1996).

\noindent
[EHZ03] P. Ebenfelt, X. Huang and D. Zaitsev,
The equivalence problem and rigidity for hypersurfaces embedded into
hyperquadrics, {\it Amer. Jour. Math.} 127 (1), 169-192, 2005.

\noindent
[Fa82] J. Faran, Maps from the two ball to the three ball, Invent. Math.
68, 441-475 (1982).

\noindent
[Fa86] J. Faran, On the linearity of proper maps between balls in the lower
dimensional case, Jour. Diff. Geom. 24, 15-17 (1986).

\noindent
[Fo89] F. Forstneric, Extending proper holomorphic mappings
of positive codimension, Invent. Math., 95, 31-62 (1989).

\noindent [Fo92] F. Forstneric, A survey on proper
holomorphic mappings, Proceeding of Year in SCVs at
Mittag-Leffler Institute, Math. Notes 38 (1992),  Princeton
University Press,  Princeton, N.J.

%\noindent [Ha91] H. Hamada, Monomial proper maps between balls of different
%dimensions,
%Bull. Kyushu Kyoritsu Univ. Fac. Engineer, 15, 41-43, (1991).

\noindent [Ha05] H. Hamada, Rational proper holomorphic maps from ${\bf{B}}^n$
into
${\bf{B}}^{2n}$,  {\it Math. Ann.}  331 (no.3), 693--711, 2005.

\noindent
[Hu99] X. Huang, On a linearity problem of proper holomorphic mappings
between balls in complex spaces of different dimensions,  Jour. of
Diff. Geom. Vol (51) No. 1, 13-33 (1999).

\noindent
[Hu01] X. Huang,  On some problems in several complex variables and CR
geometry.
First International Congress of Chinese
Mathematicians  (Beijing, 1998), 383--396, AMS/IP Stud. Adv. Math., 20,
Amer. Math. Soc., Providence, RI, 2001.

\noindent
[Hu03] X. Huang, On a Semi-Rigidity Property for Holomorphic Maps,
Asian J. Math. Vol(7) No. 4, 463-492.

\noindent [HJ01] X. Huang and S. Ji, Mapping ${\bf{B}}^n$ into ${\BB}^{2n-1}$,
Invent. Math. 145, 219-250(2001).

\noindent [HJX05] X. Huang, S. Ji, and D. Xu,
Several Results for Holomorphic Mappings from ${\bf B}^n$ into ${\BB}^N$,  {\it
Contemporary Math} 368,
AMS, 267-293 (2005). (A special issue dedicated to  Francois Treves).

\noindent
[JX04] S. Ji and D. Xu, Rational maps between $\bf{B}^n$ and $\bf{B}^N$
with geometric rank $\kappa_0\le n-2$
and minimal target dimension, Asian J. Math. Vol(8) No. 2, 233-258, 2004.

\noindent
[Mir03]  N. Mir,  Analytic regularity of CR maps into spheres.
Math. Res. Lett.  10 (no.4), 447--457, 2003.

%\noindent
%$[T62] N. Tanaka, On the pseudo-conformal geometry of hypersurfaces of
 %the space of n
%complex variables, J. Math. Soc. Japan 14, 397-429(1962).

\noindent
[W79] S. Webster, On mapping an (n+1)-ball in the complex
space, Pac. J. Math. 81, 267-272 (1979).

\bigskip
\bigskip

X. Huang and S. Ji, School of Mathematics, Wuhan University, Wuhan, Hubei
430072, China.

X. Huang, huangx@math.rutgers.edu, Department of Mathematics, Rutgers
University, New Brunswick,
NJ 08903, USA;

S. Ji,\ shanyuji@math.uh.edu,\ D. Xu,\ dekangxu@math.uh.edu,\
Department of Mathematics, University of Houston, Houston, TX, 77204.

\end{document}